\documentclass{article}

\usepackage{PRIMEarxiv}

\usepackage[utf8]{inputenc} 
\usepackage[T1]{fontenc}    
\usepackage{hyperref}       
\usepackage{url}            
\usepackage{booktabs}       
\usepackage{amsfonts}       
\usepackage{nicefrac}       
\usepackage{microtype}      
\usepackage{lipsum}
\usepackage{fancyhdr}       
\usepackage{graphicx}       
\graphicspath{{media/}}     

\title{Take-away Impartial Combinatorial Games \\on Hypergraphs and other related 
\\Geometric and Discrete Structures
}

\author{
  T. H. Molena \\
  North Carolina State University \\
  Raleigh, NC, USA\\
  \texttt{thnguy22@ncsu.edu} \\
}

\begin{document}
\maketitle

\begin{abstract}

In a Take-Away Game on hypergraphs, two players take turns to remove the vertices and the hyperedges of the hypergraphs. In each turn, a player must remove either a single vertex or a hyperedge. When a player chooses to remove one vertex, all of the hyperedges that contain the chosen vertex are also removed. When a player chooses to remove one hyperedge, only that chosen hyperedge is removed. Whoever removes the last vertex wins the game. Following from the winning strategy for the Take-Away Impartial Combinatorial Games on only Oddly Uniform or only Evenly Uniform Hypergraphs, this paper is about the new winning strategy for Take-Away Games on neither Oddly nor Evenly Uniform Hypergraphs. These neither Oddly nor Evenly Uniform Hypergraphs, however, have to satisfy the specific given requirements.
\end{abstract}


\section{The special requirements.}
In this paper, only those hypergraphs satisfying these requirements are considered.

\begin{itemize}
    
    \item In each hypergraph, these following special requirements must be satisfied:
    
    \begin{itemize}
    
        \item Each vertex must appear in at least one hyperedge of at least 2 vertices.
    
        \item All hyperedges must be under exactly one of two categories:
     
    \begin{itemize}
    
    \item Category X contains exactly one hyperedge of an even number of vertices, and all vertices in that hyperedge.
\item Category Y contains any positive number of hyperedges of exactly 3 vertices, and all vertices in those hyperedges.
    
    \end{itemize}
    
    \begin{figure}[ht]
    \centering
    \includegraphics[width=16cm]{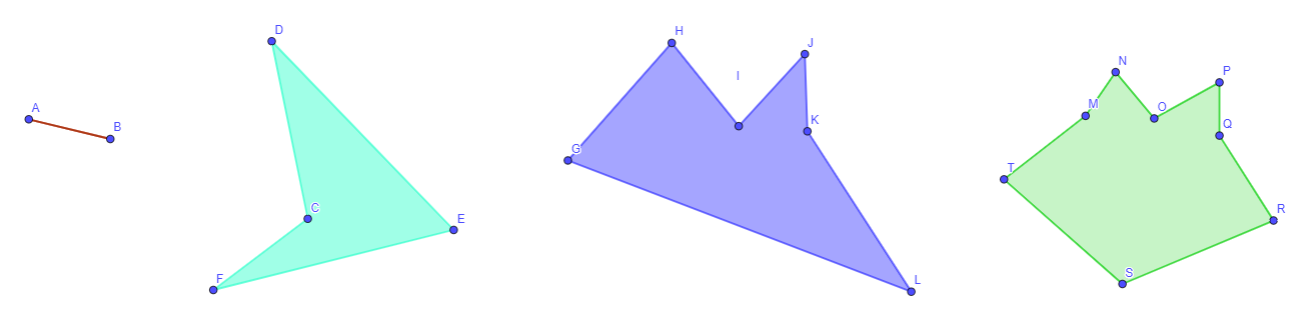}
    \caption{    Several examples of category X:
    \\Category X contains hyperedge $\{A, B\}$, vertex A, and vertex B. 
        \\Category X contains hyperedge $\{C, D, E, F\}$, vertex C, vertex D, vertex E, and vertex F.
        \\ Category X contains hyperedge \{G, H, I, J, K, L\}, vertex G, vertex H, vertex I, vertex J, vertex K, and vertex L.
        \\ Category X contains hyperedge \{M, N, O, P, Q, R, S, T\}, vertex M, vertex N, vertex O, vertex P, vertex Q, vertex R, vertex S, and vertex T.
    }
    \label{fig:my_label}
\end{figure}

    \begin{figure}[ht]
   
    \centering
    \includegraphics[width=16cm]{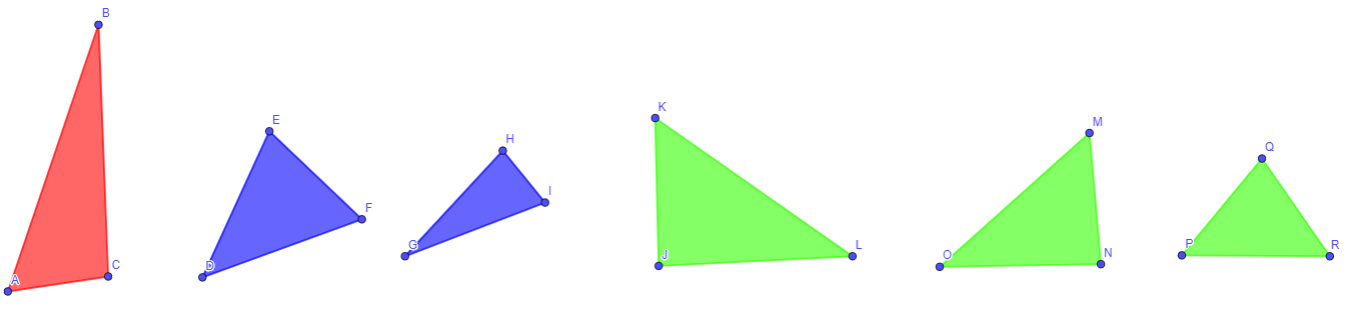}
    \caption{ 
    Several examples of category Y:
    \\Category Y contains hyperedge \{A, B, C\}, vertex A, vertex B, and vertex C.
      \\ Category Y contains hyperedge \{D, E, F\}, hyperedge \{G, H, I\}, vertex D, vertex E, vertex F, vertex G, vertex H, and vertex I.
        \\Category Y contains hyperedge \{J, K, L\}, hyperedge \{M, N, O\}, and hyperedge \{P, Q, R\}, vertex J, vertex K, vertex L, vertex M, vertex N, vertex O, vertex P, vertex Q, and vertex R.  } 
\end{figure}

        Notation:
    \begin{itemize}
        
        \item E(CatX) denotes the set of hyperedges in category X, and $|E(CatX)|$ is the cardinality of E(CatX).
    
    \item V(CatX) denotes the set of vertices in category X, and $|V(CatX)|$ is the cardinality of V(CatX). Note that V(CatX) is also the set of all vertices of an only hyperedge in E(CatX). 
        
    \item E(CatY) denotes the set of hyperedges in category Y, and $|E(CatY)|$ is the cardinality of E(CatY).
    
    \item V(CatY) denotes the set of vertices in category Y, and $|V(CatY)|$ is the cardinality of V(CatY). Note that V(CatY) is also the set of vertices of all hyperedges in E(CatY).
        
    \end{itemize}
    
    With the previous special requirements:

    \begin{itemize}

        \item $|V(CatX)|$ is even.
        \item $|E(CatX)|$ = 1.
        \item $|V(CatY)| > $ 0. 
        \item $|E(CatY)| > $ 0.
   
   \end{itemize}

    \end{itemize}
    
    \item Consider $|E(CatX)|$ and $|E(CatY)|$:

\begin{itemize}
   
   \item Any hypergraph with $|E(CatX)|$ = 1 and $|E(CatY)| > $ 0, by definition, is neither an evenly uniform hypergraph with a marked coloring nor an oddly uniform hypergraph.
   \\Therefore, any hypergraph with $|E(CatX)|$ = 1 and $|E(CatY)| > $ 0 must also satisfy these additional special requirements:

\begin{itemize}
    
    \item Consider E(CatX):
    \\Each vertex of an only hyperedge in E(CatX) must appear in 0 hyperedge, or 1 hyperedge, or 2 hyperedges in E(CatY). 
    
    \item Consider E(CatY):
    
    \begin{itemize}
      
         \item Exactly 2 out of 3 vertices in each hyperedge in E(CatY) must appear in an only hyperedge in E(CatX).
         \item The remaining 1 out of 3 vertices in each hyperedge in E(CatY) must appear in all hyperedges in E(CatY) and 0 hyperedge in E(CatX). That exactly 1 special vertex in V(CatY) is called vertex S. 
         \\In summary, vertex S must appear in all hyperedges in E(CatY) and 0 hyperedge in E(CatX). 
      
    \end{itemize}

Following from the previous additional special requirements: 

\item Consider V(CatX): 
\\Each vertex in V(CatX) must be under exactly 1 of 3 subcategories:

\begin{itemize}
 
 \item Subcategory A contains all vertices that appear in a hyperedge in E(CatX) and 1 hyperedge in E(CatY).

\item Subcategory B contains all vertices that appear in a hyperedge in E(CatX) and 2 hyperedges in E(CatY).

\item Subcategory C contains all vertices that appear in a hyperedge in E(CatX) and 0 hyperedge in E(CatY).

\end{itemize}

\item Consider V(CatY): 
\\Following from the previous parts, V(CatY) only contains a special vertex (vertex S), all vertices in subcategory A (if any), and all vertices in subcategory B (if any).

\end{itemize}
\end{itemize}
\end{itemize}

\begin{figure}[h]
    \centering
    \includegraphics[width=16cm]{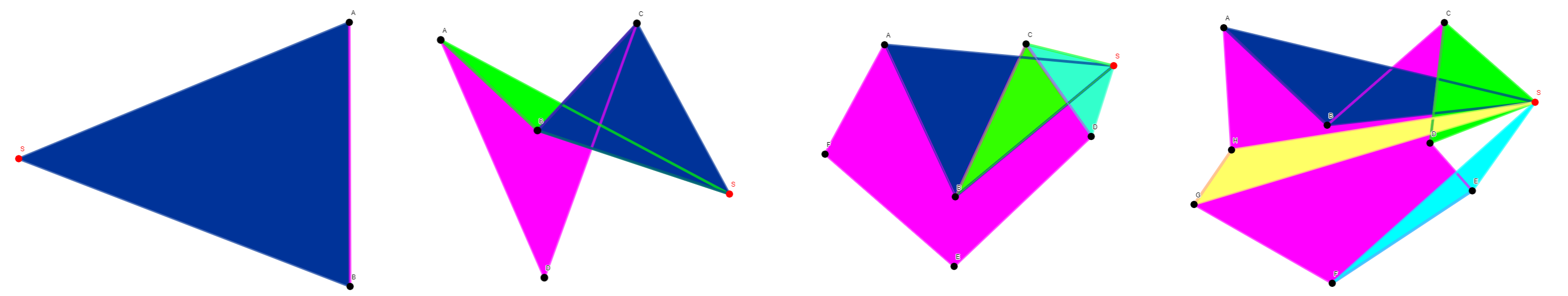}
    \caption{ 
    Several hypergraphs with $|E(CatX)|$ = 1 and $|E(CatY)| >$ 0:
    \\Category X contains a pink hyperedge, and all black vertices in each hypergraph.
    \\Category Y contains all hyperedges not in category X (not in pink), all vertices in those hyperedges, and a red vertex S in each hypergraph.
    \\The first figure (from left to right): 
          \\Category X contains hyperedge \{A, B\}, and all vertices in that hyperedge.
       \\Category Y contains hyperedge \{S, A, B\}, and all vertices in that hyperedge.
   \\The second figure (from left to right):
        \\Category X contains hyperedge \{A, B, C, D\}, and all vertices in that hyperedge
\\Category Y contains hyperedge \{S, A, B\}, hyperedge \{S, B, C\}, and all vertices in those hyperedges.
\\The third figure (from left to right):
    \\Category X contains hyperedge\{A, B, C, D, E, F\}, and all vertices in that hyperedge.
    \\Category Y contains hyperedge \{S, A, B\}, hyperedge \{S, B, C\}, hyperedge \{S, C, D\}, and all vertices in those hyperedges.
    \\The fourth figure (from left to right):
    \\Category X contains hyperedge\{A, B, C, D, E, F, G, H\}, and all vertices in that hyperedge.
    \\Category Y contains hyperedge \{S, A, B\}, hyperedge \{S, C, D\}, hyperedge \{S, E, F\}, hyperedge \{S, G, H\}, and all vertices in those hyperedges.  
    }
\end{figure}

For example: In a third hypergraph (from left to right) in Figure 3:
\\Consider category X: 
\begin{itemize}
    \item E(CatX) contains hyperedge\{A, B, C, D, E, F\}
    \item V(CatX) contains vertex A, vertex B, vertex C, vertex D, vertex E, and vertex F.
\end{itemize}
Consider category Y:
\begin{itemize}
    \item E(CatY) contains hyperedge \{S, A, B\}, hyperedge \{S, B, C\}, hyperedge \{S, C, D\} 
    \item V(CatY) contains all vertices in the hyperedges in E(CatY).
\end{itemize}
Since vertex A and vertex D appear in an only hyperedge in E(CatX) and 1 hyperedge in E(CatY), vertex A and vertex D are in subcategory A.
\\Since vertex B and vertex C appear in an only hyperedge in E(CatX) and 2 hyperedges in E(CatY), vertex B and vertex C are in subcategory B.
\\Since vertex E and vertex F appear in an only hyperedge in E(CatX) and 0 hyperedge in E(CatY), vertex E and vertex F are in subcategory C.

\vspace{0.2cm}

For hypergraphs with 
\begin{itemize}
    \item Either an only hyperedge (with all vertices in that hyperedge) in category X, and 0 hyperedge in category Y.
    \item Or only hyperedges (with all vertices in those hyperedges) in category Y, and 0 hyperedge in category X. 
\end{itemize}
the previous additional special requirements do not apply. Therefore, in those cases, there does not exist a special vertex S nor any subcategory. Instead, the winning strategy was already found in those cases. By definition: 

\begin{itemize}
    
    \item Any hypergraph with an only hyperedge (with all vertices in that hyperedge) in category X, and 0 hyperedge in category Y, is an evenly uniform hypergraph with a marked coloring.
    
    \item Any hypergraph with only hyperedges (with all vertices in those hyperedges) in category Y, and 0 hyperedge in category X is an oddly uniform hypergraph.
    
\end{itemize}

Take-Away impartial combinatorial games on both evenly uniform hypergraphs with a marked coloring, and oddly uniform hypergraphs were already solved \cite{barnardkristen}. Therefore, in this article, the focus is only on the Take-Away impartial combinatorial games on any hypergraphs with $|E(CatX)|$ = 1 and $|E(CatY)| > $ 0 that satisfies the special requirements in Part 1.

\section{Lemmas, theorems, and corollary.}
\label{sec:headings}

\subsection{Lemma 1}{}
There does not exist a hypergraph where all vertices in V(CatX) are in subcategory C.

\begin{paragraph}{Proof}
\hspace{0cm}
\\To the contrary, assume that there exists a hypergraph where all vertices in V(CatX) are in subcategory C.
\\Vertices in subcategory C appear in a hyperedge in E(CatX) and 0 hyperedge in E(CatY). Furthermore, from the assumption and the additional special requirements, V(CatY) only contains a vertex S. 
\\It means that when all vertices in V(CatX) are in subcategory C, a hypergraph contains a hyperedge (and all vertices in that hyperedge) in category X. 
\\Consequently, the previous additional special requirements do not apply. Therefore, in those cases, there does not exist a special vertex S nor any subcategory. 
\\Thus, by contradiction, there does not exist a hypergraph where all vertices in V(CatX) are in subcategory C. 

\end{paragraph}


\subsection{Lemma 2}{}
Prove that $|E(CatY)|$ = $|V(CatX)|$ for all hypergraphs where all vertices in V(CatX) are in subcategory B. 

\begin{paragraph}{Proof}
\hspace{0cm}
\\Consider all hypergraphs where all vertices in V(CatX) are in subcategory B. 
\\Vertices in subcategory B appear in a hyperedge in E(CatX) and 2 hyperedges in E(CatY). Furthermore, from the hypothesis and the additional special requirements, V(CatY) only contains a vertex S, and all vertices in subcategory B. Suppose that 2 out of 3 vertices in each hyperedge in E(CatY) can be the same vertices, it leads to $|E(CatY)|$ = $|V(CatX)|$. With the fact that 2 out of 3 vertices in each hyperedge in E(CatY) must be distinct, with a vertex S and the same number of vertices in subcategory B , there are a lot of different combinations to make valid hyperedges in E(CatY), so $|E(CatY)|$ and $|V(CatX)|$ both remain unchanged, and $|E(CatY)|$ = $|V(CatX)|$.
\\Thus, $|E(CatY)|$ = $|V(CatX)|$ for all hypergraphs where all vertices in V(CatX) are in subcategory B.
\\There is an additional remark that with a given $|V(CatX)|$, as exactly 2 out of 3 vertices in each hyperedge in E(CatY) must appear in a hyperedge in E(CatX), the maximum possible number that $|E(CatY)|$ can reach is the number which is equal to $|V(CatX)|$. Consequently, for all hypergraphs where all vertices in V(CatX) are in subcategory B, $|E(CatY)|$ reaches the maximum possible number with the given $|V(CatX)|$, which is equal to $|V(CatX)|$ .   
\end{paragraph}

\begin{figure}[ht]
    \centering
    \includegraphics[width=15cm]{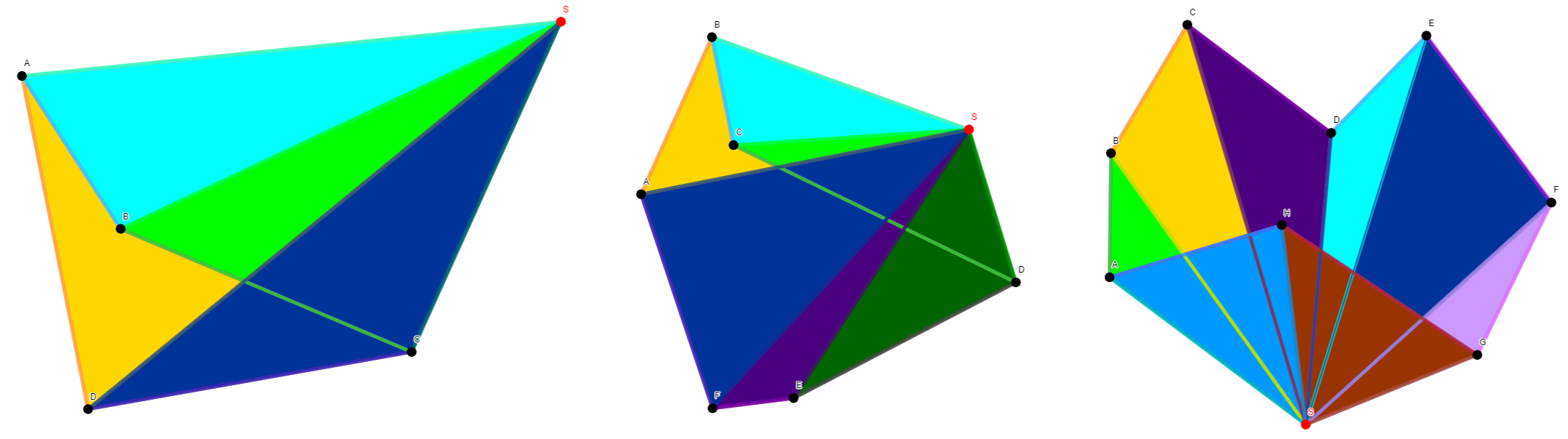}
    \caption{ 
    Several hypergraphs with $|E(CatY)|$ = $|V(CatX)|$.
    }
\end{figure}

\subsection{Lemma 3}{}
Prove that $|E(CatY)| \geq$ 4 for all hypergraphs where all vertices in V(CatX) are in subcategory B.

\begin{paragraph}{Proof}
\hspace{0cm}
\\To the contrary, assume that 0 $ \leq |E(CatY)| < $ 4 in hypergraphs where all vertices in V(CatX) are in subcategory B. 
\\From Lemma 2, $|E(CatY)|$ = $|V(CatX)|$ for all hypergraphs where all vertices in V(CatX) are in subcategory B. From the special requirement, since $|V(CatX)|$ is even, $|E(CatY)|$ is even.
\\Since $|E(CatY)|$ is even and 0 $ \leq |E(CatY)| < $ 4, there are only 2 cases to consider: either $|E(CatY)|$ = 0 or $|E(CatY)|$ = 2. 

\begin{itemize}

    \item Case 1: $|E(CatY)|$ = 0. 
    \\For hyperedges with only a hyperedge (and all vertices in that hyperedge) in cateogry X, the previous additional special requirements do not apply. Therefore, in those cases, there does not exist a special vertex S nor any subcategory.
    \\Consequently, this case reaches a contradiction.
    
    \item Case 2: $|E(CatY)|$ = 2. 
    \\From the same argument in Lemma 2, for all hypergraphs where all vertices in V(CatX) are in subcategory B, V(CatY) only contains a vertex S, and all vertices in subcategory B. Vertices in subcategory B appears in a hyperedge in E(CatX) and 2 hyperedges in E(CatY). Consequently, in case 2, those 2 hyperedges of 3 distinct vertices in E(CatY) must be exactly the same. Therefore, there is only 1 valid hyperedge of 3 vertices in E(CatY), which reaches a contradiction.

\end{itemize}
In conclusion, a contradiction is found in either case. 
\\Thus, by contradiction, $|E(CatY)| \geq$ 4 for all hypergraphs where all vertices in V(CatX) are in subcategory B.

\end{paragraph}

\begin{figure}[ht]
    \centering
    \includegraphics[width=15cm]{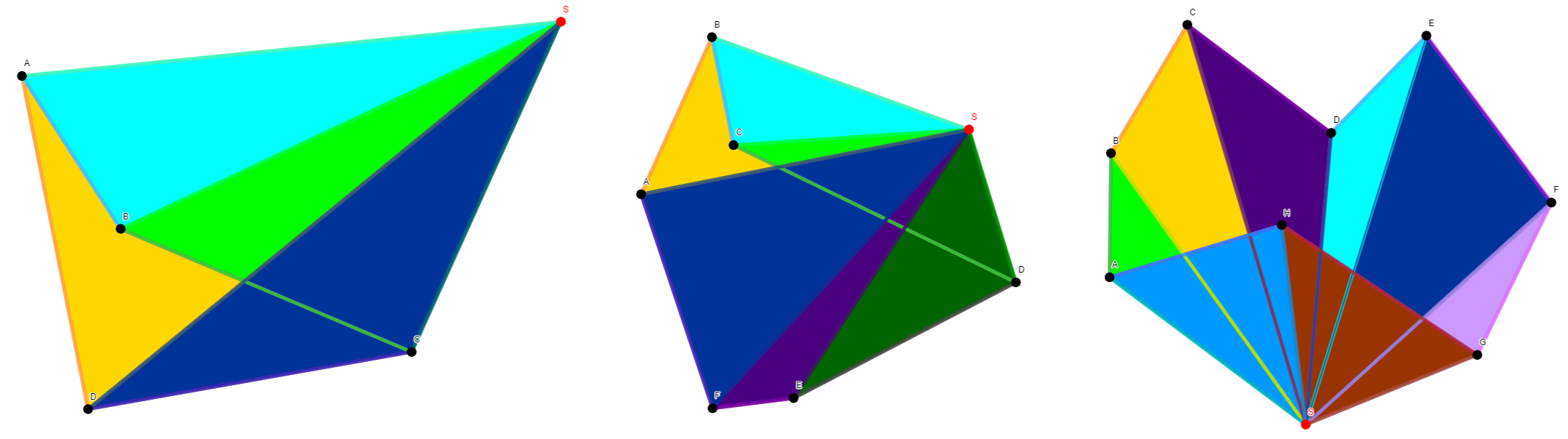}
    \caption{ 
Several hypergraphs with all vertices in $|V(CatX)|$ are in subcategory B.    }
\end{figure}

\subsection{Lemma 4}{}
There does not exist a hypergraph where all vertices in V(CatX) are in subcategory B and subcategory C.

\begin{paragraph}{Proof}

\hspace{0cm}
\\To the contrary, assume that there exists a hypergraph where all vertices in V(CatX) are in subcategory B and subcategory C.
\\Vertices in subcategory C appear in a hyperedge in E(CatX)
and 0 hyperedge in E(CatY). Vertices in subcategory B appear in a hyperedge in E(CatX) and 2 hyperedges in E(CatY)
Furthermore, from the hypothesis and the additional special requirements, V(CatY) only contains a special vertex
(vertex S), and all vertices in subcategory B.
\\From Lemma 2,  $|E(CatY)|$ = $|V(CatX)|$ for all hypergraphs where V(CatY) only contains a vertex S, and all vertices in subcategory B. 
\\As there is 0 vertex that appears in a hyperedge in E(CatX) and 0 hyperedge in E(CatY), there is 0 vertex in subcategory C.
\\Thus, by contradiction, there does not exist a hypergraph where all vertices in V(CatX) are in subcategory B and subcategory C.

\end{paragraph}

\subsection{Lemma 5}{}
Prove that $|E(CatY)| \geq$ 3 for all hypergraphs where all vertices in V(CatX) are in subcategory A and subcategory B.

\begin{paragraph}{Proof}
\hspace{0cm}
\\To the contrary, assume that 0 $ \leq |E(CatY)| < $ 3 in hypergraphs where all vertices in V(CatX) are in subcategory A and subcategory B. 
\begin{itemize}
    \item Case 1: $|E(CatY)|$ = 0. 
    \\For hypergraphs with only a hyperedge (and all vertices in that hyperedge) in category X, the previous additional special requirements do not apply. Therefore,
    in those cases, there does not exist a special vertex S nor any subcategory.
    \\Consequently, this case reaches a contradiction.
    \item Case 2: $|E(CatY)|$ = 1.
    \\From the hypothesis, and the additional special requirements, V(CatY) must contain only a vertex S and 2 vertices in subcategory A.
    \\Therefore, V(CatX) can only contain all vertices in subcategory A and all vertices in category C (if any).
    \\Consequently, this case reaches a contradiction.
    \item Case 3: $|E(CatY)|$ = 2. 
    \\From the hypothesis, and the additional special requirements, V(CatY) can only contain a vertex S, 1 vertex in subcategory B and 2 vertices in subcategory A. However, this leads to V(CatX) = 3 (2 vertices in subcategory A and 1 vertex in subcategory B), which is not an even number of vertices.
    \\Consequently, this case reaches a contradiction.
\end{itemize}
In conclusion, a contradiction is found in all 3 cases.
\\Thus, by contradiction, $|E(CatY)| \geq$ 3 for all hypergraphs where all vertices in V(CatX) are in subcategory A and subcategory B.
\end{paragraph}

\subsection{Lemma 6}{}
From Lemma 1 and Lemma 4, there is a corollary:
\\With 3 subcategories (subcategory A, subcategory B, and subcategory C), there are exactly 5 distinct groups (counted in Roman numerals for easier use later) for all hypergraphs of all vertices in V(CatX):
\begin{itemize}
    \item Group I: Hypergraphs where all vertices in V(CatX) are in subcategory A.
    \item Group II: Hypergraphs where all vertices in V(CatX) are in subcategory B. For all hypergraphs in this case, from Lemma 2 and Lemma 3, $|E(CatY)|$ = $|V(CatX)|$, and $|E(CatY)|$ is even (as $|V(CatX)|$ is even).   
    \item Group III: Hypergraphs where all vertices in V(CatX) are in subcategory A and subcategory B.
    \item Group IV: Hypergraphs where all vertices in V(CatX) are in subcategory A and subcategory C.
    \item Group V: Hypergraphs where all vertices in V(CatX) are in subcategory A, subcategory B, and subcategory C.
\end{itemize}
From Corollary 6, there is a remark:
\begin{itemize}
    \item Except all hypergraphs in group II, for any hypergraph in 1 of 4 remaining groups, there is always at least a vertex in V(CatX) that is in subcategory A.
    \item Consider all hypergraphs in group III, group IV, and group V. For any hypergraph in 1 of those 3 groups, there is always at least a vertex in V(CatX) that is in subcategory B or there is always at least a vertex in V(CatX) that is in subcategory C.
\end{itemize}

\begin{paragraph}{Theorem 7}
\hspace{1cm}
\\If T is a hypergraph where $|E(CatX)|$ = 1 and $|E(CatY)|$ is odd, then g(T) = 1. 
\end{paragraph}

\begin{paragraph}{Proof}
\hspace{1cm}

Proof by induction.
\\Base case: Considering a hypergraph T where $|E(CatX)|$ = 1 and $|E(CatY)|$ = 1. 
\begin{itemize}
    \item Since $|V(CatX)|$ is even, and there is a vertex S (not in category X but in category Y). Consequently, the vertex set $|V(T)|$ is odd. 
    \item Since $|E(CatY)|$ = 1, and by the additional special requirements, category Y contains a vertex S and 2 vertices in subcategory A. 
\end{itemize}
Consequently, all hypergraphs in group I and group IV are valid with this hypothesis. 
\\Therefore, there is always at least a vertex in subcategory A in all valid hypergraphs.
\\By the rule of the game, either 1 vertex or 1 hyperedge is removed at a time.
\begin{itemize}
    \item Case 1: Remove a vertex S. 
    \\After a vertex S is removed:
    
    \begin{itemize}
        \item $|V(CatX)|$ is even.
        \item $|E(CatX)|$ = 1, and $|E(CatY)|$ = 0.
    \end{itemize}
    
    By definition, the resulting hypergraph $T_{1}$ is an evenly uniform hypergraph with a marked coloring with: 
    
    \begin{itemize}
        \item The resulting vertex set $|V(T_{1})|$ is even.
        \item The resulting hyperedge set $|E(T_{1})|$ is odd. 
    \end{itemize}
    By Theorem 2.8, g($T_{1}$) = 2.

    \item Case 2: Remove a vertex in subcategory A.
    \\After a vertex in subcategory A is removed: 
    
    \begin{itemize}
        \item $|V(CatX)|$ is odd, and there is a vertex S (not in category X but in category Y).
        \item $|E(CatX)|$ = 0, and $|E(CatY)|$ = 0.
    \end{itemize}
    
    By definition, the resulting hypergraph $T_{2}$ is an oddly uniform hypergraph with: 
    
    \begin{itemize}
        \item The resulting vertex set $|V(T_{2})|$ is even.
        \item The resulting hyperedge set $E(T_{2})$ is even.
    \end{itemize}
    By Theorem 2.4, g($T_{2}$) = 0. 
    
    \item Case 3: Remove a vertex in subcategory C (if any). 
    \\After a vertex in subcategory C is removed:
    
    \begin{itemize}
        \item $|V(CatX)|$ is odd, and there is a vertex S (not in category X but in category Y).
        \item $|E(CatX)|$ = 0, and $|E(CatY)|$ = 1.
    \end{itemize}

    By definition, the resulting hypergraph $T_{3}$ is an oddly uniform hypergraph with: 
    \begin{itemize}
        \item The resulting vertex set $|V(T_{3})|$ is even. 
        \item The resulting hyperedge set $|E(T_{3})|$ is odd. 
    \end{itemize}
    By Theorem 2.4, g($T_{3}$) = 3. 
    
    \item Case 4: Remove a hyperedge in E(CatX).
    \\After a hyperedge in E(CatX) is removed, 
    
    \begin{itemize}
        \item $|V(CatX)|$ is even, and there is a vertex S (not in category X but in category Y).
        \item $|E(CatX)|$ = 0, and $|E(CatY)|$ = 1.
    \end{itemize}
    
    By definition, the resulting hypergraph $T_{4}$ is an oddly uniform hypergraph with:
    \begin{itemize}
        \item The resulting vertex set $|V(T_{4}|$ is odd.
        \item The resulting hyperedge $|E(T_{4})|$ is odd. 
    \end{itemize} 
    By Theorem 2.4, g($T_{4}$) = 2. 
    
    \item Case 5: Remove a hyperedge in E(CatY).
    After a hyperedge in E(CatY) is removed: 
    
    \begin{itemize}
        \item $|V(CatX)|$ is even, and there is a vertex S (not in category X but in category Y).
        \item $|E(CatX)|$ = 1, and $|E(CatY)|$ = 0.
    \end{itemize}
    
    By definition, the resulting hypergraph $T_{5}$ is an evenly uniform hypergraph with a marked coloring with: 
    \begin{itemize}
        \item The resulting vertex set $|V(T_{5}|$ is odd.
        \item The resulting hyperedge set $|E(T_{5})|$ is odd.  
    \end{itemize}
    By Theorem 2.8, g($T_{5}$) = 3.
\end{itemize}

\begin{itemize}
    \item Consider all hypergraphs T in group I. By Sprague–Grundy Theorem, g(T) = mex\{2,0,2,3\} = 1. 
    
    \item Consider all hypergraphs T in group IV. By Sprague–Grundy Theorem, g(T) = mex\{2,0,3,2,3\} = 1.
    
\end{itemize}
Therefore, for all hypergraphs T in either group I or group IV, g(T) = 1.

\vspace{0.1in}
Inductive hypothesis: assume that if T is a hypergraph where $|E(CatX)|$ = 1 and $|E(CatY)|$ = n (n is odd), then g(T) = 1.
\\Consider a hypergraph T where $|E(CatX)|$ = 1 and $|E(CatY)|$ = n + 1 (n + 1 is odd). 
\\Prove that $g(T)$ = 1.

\begin{itemize}
    \item Since $|V(CatX)|$ is even, and there is a vertex S (not in category X but in category Y). Consequently, the vertex set $|V(T)|$ is odd.
    \item The hyperedge set E(T) contains $|E(CatX)|$ = 1 and $|E(CatY)|$ is odd.
\end{itemize}

By Corollary 6, all hypergraphs in group II are not valid in this Theorem (as the inductive hypothesis states that $|E(CatY)|$ is odd). Only all hypergraphs in group I, group III, group IV, and group V are valid.
\\Therefore, by Corollary 6, there is always at least a vertex in subcategory A in all valid hypergraphs.
\\By the rule of the game, either 1 vertex or 1 hyperedge is removed at a time.
\begin{itemize}
    \item Case 1: Remove a vertex S.
    \\After a vertex S is removed: 
    \begin{itemize}
        \item $|V(CatX)|$ is even.
        \item $|E(CatX)|$ = 1, and $|E(CatY)|$ = 0.
    \end{itemize}
    By definition, the resulting hypergraph $T_{1}$ is an evenly uniform hypergraph with a marked coloring with:
    \begin{itemize}
        \item The resulting vertex set $|V(T_{1})|$ is even.
        \item The resulting hyperedge set $|E(T_{1})|$ is odd.  
    \end{itemize}
    By Theorem 2.8, g($T_{1}$) = 2.

    \item Case 2: Remove a vertex in subcategory A. 
    \\After a vertex in subcategory A is removed: 
    \begin{itemize}
        \item $|V(CatX)|$ is odd, and there is a vertex S (not in category X but in category Y).
        \item $|E(CatX)|$ = 0, and $|E(CatY)|$ is even.
    \end{itemize}
    By definition, the resulting hypergraph $T_{2}$ is an oddly uniform hypergraph with:
    \begin{itemize}
        \item The resulting vertex set $|V(T_{2})|$ is even.
        \item The resulting hyperedge set $|E(T_{2})|$ is even.
    \end{itemize}
    By Theorem 2.4, g($T_{2}$) = 0.

    \item Case 3: Remove a vertex in subcategory B (if any). 
    \\After a vertex in subcategory B is removed:
    \begin{itemize}
    \item $|V(CatX)|$ is odd, and there is a vertex S (not in category X but in category Y).
    \item $|E(CatX)|$ = 0, and $|E(CatY)|$ is odd.
    \end{itemize}
    By definition, the resulting hypergraph $T_{3}$ is an oddly uniform hypergraph with:
    \begin{itemize}
        \item The resulting vertex set $|V(T_{3})|$ is even.
        \item The resulting hyperedge set $|E(T_{3})|$ is odd. 
    \end{itemize}
    By Theorem 2.4, g($T_{3}$) = 3.
    \item Case 4: Remove a vertex in subcategory C (if any). 
    \\After a vertex in subcategory C is removed:
    \begin{itemize}
    \item $|V(CatX)|$ is odd, and there is a vertex S (not in category X but in category Y).
    \item $|E(CatX)|$ = 0, and $|E(CatY)|$ is odd.
    \end{itemize}
    
    By definition, the resulting hypergraph $T_{4}$ is an oddly uniform hypergraph with: 
    
    \begin{itemize}
        \item The resulting vertex set $|V(T_{4})|$ is even.
        \item The resulting hyperedge set $|E(T_{4})|$ is odd. 
    \end{itemize}
    By Theorem 2.4, g($T_{4}$) = 3.
    
    \item Case 5: Remove a hyperedge in E(CatX).
    \\After a hyperedge in E(CatX) is removed:
    
    \begin{itemize}
    \item $|V(CatX)|$ is even, and there is a vertex S (not in category X but in category Y).
    \item $|E(CatX)|$ = 0, and $|E(CatY)|$ is odd.
    \end{itemize}
    By definition, the resulting hypergraph $T_{5}$ is an oddly uniform hypergraph with:
    
    \begin{itemize}
        \item The resulting vertex set $|V(T_{5}|$ is odd.
        \item The resulting hyperedge set $|E(T_{5})|$ is odd. 
    \end{itemize}
    By Theorem 2.4, g($T_{5}$) = 2.
    
    \item Case 6: Remove a hyperedge in E(CatY).
    \\After a hyperedge in E(CatY) is removed:
    
    \begin{itemize}
    \item Since $|V(CatX)|$ is even, and there is a vertex S (not in category X but in category Y), $|V(T_{6}|$ is odd.
    \item $|E(CatX)|$ = 1, and $|E(CatY)|$ is even.
    \\By definition, the resulting hypergraph $T_{6}$ is neither an oddly uniform hypergraph nor an evenly uniform hypergraph with a marked coloring. 
    \\Starting from here, $|E(CatY)|$ is even, but from the starting note above, all hypergraphs in group II are already not valid. Therefore, only all hypergraphs in group I, group III, group IV, and group V are valid.
\\Therefore, by Corollary 6, there is always at least a vertex in subcategory A in all valid hypergraphs.
\\Continue with the rule of the game, either 1 vertex or 1 hyperedge is removed at a time.

    \begin{itemize}
    \item Subcase 1: Remove a vertex S.
    \\After a vertex S is removed: 
    \begin{itemize}
        \item $|V(CatX)|$ is even.
        \item $|E(CatX)|$ = 1, and $|E(CatY)|$ = 0.
    \end{itemize}
    By definition, the resulting hypergraph $T_{6_1}$ is an evenly uniform hypergraph with a marked coloring with:
    \begin{itemize}
        \item The resulting vertex set $|V(T_{6_1})|$ is even.
        \item The resulting hyperedge set $|E(T_{6_1})|$ is odd.  
    \end{itemize}
    By Theorem 2.8, g($T_{6_1}$) = 2.
    
   \item Subcase 2: Remove a vertex in subcategory A. 
    \\After a vertex in subcategory A is removed: 
    \begin{itemize}
        \item $|V(CatX)|$ is odd, and there is a vertex S (not in category X but in category Y).
        \item $|E(CatX)|$ = 0, and $|E(CatY)|$ is odd.
    \end{itemize}
    By definition, the resulting hypergraph $T_{6_2}$ is an oddly uniform hypergraph with:
    \begin{itemize}
        \item The resulting vertex set $|V(T_{6_2})|$ is even.
        \item The resulting hyperedge set $|E(T_{6_2})|$ is odd.
    \end{itemize}
    By Theorem 2.4, g($T_{6_2}$) = 3. 
  
    \item Subcase 3: Remove a vertex in subcategory B (if any). 
    \\After a vertex in subcategory B is removed: 
    \begin{itemize}
        \item $|V(CatX)|$ is odd, and there is a vertex S (not in category X but in category Y).
        \item $|E(CatX)|$ = 0, and $|E(CatY)|$ is even.
    \end{itemize}
    By definition, the resulting hypergraph $T_{6_3}$ is an oddly uniform hypergraph with:
    \begin{itemize}
        \item The resulting vertex set $|V(T_{6_3})|$ is even.
        \item The resulting hyperedge set $|E(T_{6_3})|$ is even.
    \end{itemize}
    By Theorem 2.4, g($T_{6_3}$) = 0.
    
    \item Subcase 4: Remove a vertex in subcategory C (if any). 
    \\After a vertex in subcategory C is removed: 
    \begin{itemize}
        \item $|V(CatX)|$ is odd, and there is a vertex S (not in category X but in category Y).
        \item $|E(CatX)|$ = 0, and $|E(CatY)|$ is even.
    \end{itemize}
    By definition, the resulting hypergraph $T_{6_4}$ is an oddly uniform hypergraph with:
    \begin{itemize}
        \item The resulting vertex set $|V(T_{6_4})|$ is even.
        \item The resulting hyperedge set $|E(T_{6_4})|$ is even.
    \end{itemize}
    By Theorem 2.4, g($T_{6_4}$) = 0. 
    
    \item Subcase 5: Remove a hyperedge in E(CatX).
    \\After a hyperedge in E(CatX) is removed:
    
    \begin{itemize}
    \item $|V(CatX)|$ is even, and there is a vertex S (not in category X but in category Y).
    \item $|E(CatX)|$ = 0, and $|E(CatY)|$ is odd.
    \end{itemize}
    By definition, the resulting hypergraph $T_{6_5}$ is an oddly uniform hypergraph with:
    
    \begin{itemize}
        \item The resulting vertex set $|V(T_{6_5}|$ is odd.
        \item The resulting hyperedge set $|E(T_{6_5})|$ is even. 
    \end{itemize}
    By Theorem 2.4, g($T_{6_5}$) = 1.
    
    \item Subcase 6: Remove a hyperedge in E(CatY).
    \\After a hyperedge in E(CatY) is removed:
    
    \begin{itemize}
    \item Since $|V(CatX)|$ is even, and there is a vertex S (not in category X but in category Y), $|V(T_{6_6})|$ is odd.
    \item $|E(CatX)|$ = 1, and $|E(CatY)|$ is odd.
    \end{itemize}
    By the Inductive Hypothesis, g($T_{6_6}$) = 1. 
\end{itemize}
\end{itemize}
\begin{itemize}
    \item Consider all hypergraphs in group I. By Sprague-Grundy Theorem, $g(T_{6})$ = mex \{2, 3, 1, 1\} = 0.
    \item Consider all hypergraphs in group III, group IV, and group V. By Sprague-Grundy Theorem, $g(T_{6})$ = mex \{2, 3, 0, 1, 1\} = 4.
\end{itemize}
Therefore, for all hypergraphs in all groups, either $g(T_{6})$ = 0 or $g(T_{6})$ = 4.
\end{itemize}

\begin{itemize}
    \item Consider all hypergraphs in group I. By Sprague-Grundy Theorem, $g(T_{6})$ = mex \{2, 0, 2, 0\} = 1, or by Sprague-Grundy Theorem, $g(T_{6})$ = mex \{2, 0, 2, 4\} = 1.
    \item Consider all hypergraphs in group III, group IV, and group V. By Sprague-Grundy Theorem, $g(T_{6})$ = mex \{2, 0, 3, 2, 0\} = 1 or by Sprague-Grundy Theorem, $g(T_{6})$ = mex \{2, 0, 3, 2, 4\} = 1.
\end{itemize}   
Therefore, for all hypergraphs in all groups, g(T) = 1. 
\\Thus, by induction, if T is a hypergraph where $|E(CatX)|$ = 1 and $|E(CatY)|$ is odd, then g(T) = 1.
\end{paragraph}

\begin{paragraph}{Theorem 8}
\hspace{0.1cm}
\\For all hypergraphs T in group III, group IV, and group V where $|E(CatX)|$ = 1 and $|E(CatY)|$ is a positive even number, then g(T) = 4. 
\end{paragraph}

\begin{paragraph}{Proof}
\hspace{1cm}
\\Proof by induction.
\\Base case: Consider a hypergraph T in group III, or group IV, or group V where $|E(CatX)|$ = 1 and $|E(CatY)|$ = 2. 
\begin{itemize}
    \item Since $|V(CatX)|$ is even, and there is a vertex S (not in category X but in category Y). Consequently, the vertex set $|V(T)|$ is odd. 
    \item By Corollary 6, since a hypergraph T is in group III, or group IV, or group V, there is always at least a vertex in subcategory A in all valid hypergraphs.
    \item Furthermore, by Corollary 6, since a hypergraph T is in group III, or group IV, or group V, there is always at least a vertex in subcategory B, or there is always at least a vertex in subcategory C in all valid hypergraphs.
\end{itemize}

By the rule of the game, either 1 vertex or 1 hyperedge is removed
at a time.
\begin{itemize}
    \item Case 1: Remove a vertex S. 
    \\After a vertex S is removed:
    
    \begin{itemize}
        \item $|V(CatX)|$ is even.
        \item $|E(CatX)|$ = 1, and $|E(CatY)|$ = 0.
    \end{itemize}
    
    By definition, the resulting hypergraph $T_{1}$ is an evenly uniform hypergraph with a marked coloring with: 
    
    \begin{itemize}
        \item The resulting vertex set $|V(T_{1})|$ is even.
        \item The resulting hyperedge set $|E(T_{1})|$ is odd. 
    \end{itemize}
    By Theorem 2.8, g($T_{1}$) = 2.

    \item Case 2: Remove a vertex in subcategory A.
    \\After a vertex in subcategory A is removed: 
    
    \begin{itemize}
        \item $|V(CatX)|$ is odd, and there is a vertex S (not in category X but in category Y).
        \item $|E(CatX)|$ = 0, and $|E(CatY)|$ = 1.
    \end{itemize}
    
    By definition, the resulting hypergraph $T_{2}$ is an oddly uniform hypergraph with: 
    
    \begin{itemize}
        \item The resulting vertex set $|V(T_{2})|$ is even.
        \item The resulting hyperedge set $E(T_{2})$ is odd.
    \end{itemize}
    By Theorem 2.4, g($T_{2}$) = 3.

    \item Case 3: Remove a vertex in subcategory B (if any).
     \\After a vertex in subcategory B is removed:
    \begin{itemize}
    \item $|V(CatX)|$ is odd, and there is a vertex S (not in category X but in category Y).
    \item $|E(CatX)|$ = 0, and $|E(CatY)|$ = 0.
    \end{itemize}
    By definition, the resulting hypergraph $T_{3}$ is an oddly uniform hypergraph with:
    \begin{itemize}
        \item The resulting vertex set $|V(T_{3})|$ is even.
        \item The resulting hyperedge set $|E(T_{3})|$ is even. 
    \end{itemize}
    By Theorem 2.4, g($T_{3}$) = 0.

    \item Case 4: Remove a vertex in subcategory C (if any). 
    \\After a vertex in subcategory C is removed:
    \begin{itemize}
    \item $|V(CatX)|$ is odd, and there is a vertex S (not in category X but in category Y).
    \item $|E(CatX)|$ = 0, and $|E(CatY)|$ = 0.
    \end{itemize}
    By definition, the resulting hypergraph $T_{4}$ is an oddly uniform hypergraph with:
    \begin{itemize}
        \item The resulting vertex set $|V(T_{4})|$ is even.
        \item The resulting hyperedge set $|E(T_{4})|$ is even. 
    \end{itemize}
    By Theorem 2.4, g($T_{4}$) = 0.
    
    \item Case 5: Remove a hyperedge in E(CatX).
    \\After a hyperedge in E(CatX) is removed:
    
    \begin{itemize}
    \item $|V(CatX)|$ is even, and there is a vertex S (not in category X but in category Y).
    \item $|E(CatX)|$ = 0, and $|E(CatY)|$ = 2.
    \end{itemize}
    By definition, the resulting hypergraph $T_{5}$ is an oddly uniform hypergraph with:
    
    \begin{itemize}
        \item The resulting vertex set $|V(T_{5}|$ is odd.
        \item The resulting hyperedge set $|E(T_{5})|$ is even. 
    \end{itemize}
    By Theorem 2.4, g($T_{5}$) = 1.

    \item Case 6: Remove a hyperedge in E(CatY).
    \\After a hyperedge in E(CatY) is removed:
    
    \begin{itemize}
    \item $|V(CatX)|$ is even, and there is a vertex S (not in category X but in category Y).
    \item $|E(CatX)|$ = 1, and $|E(CatY)|$ = 1.
    \end{itemize}
    By Theorem 7, g($T_{6}$) = 1. 
\end{itemize}
Consider all hypergraphs T in group III, group IV, and group V. By Sprague-Grundy Theorem, $g(T)$ = mex \{2, 3, 0, 1, 1\} = 4.

\vspace{0.1in}

Inductive hypothesis: Consider all hypergraphs T in group III, group IV, and group V where $|E(CatX)|$ = 1 and $|E(CatY)|$ = n (n is a positive even number), then g(T) = 4.
\\Consider all hypergraphs T in group III, group IV, and group V where $|E(CatX)|$ = 1 and $|E(CatY)|$ = n+1 (n+1 is a positive even number). Prove that g(T) = 4.
\begin{itemize}
    \item Since $|V(CatX)|$ is even, and there is a vertex S (not in category X but in category Y). Consequently, the vertex set $|V(T)|$ is odd.
    \item The hyperedge set E(T) contains $|E(CatX)|$ = 1 and $|E(CatY)|$ is a positive even number.
\end{itemize}
By Corollary 6, since a hypergraph T is group III, or group IV, or group V, there is always at least a vertex in subcategory A in all valid hypergraphs.
\\Furthermore, by corollary 6, since a hypergraph T is group III, or group IV, or group V, there is always at least a vertex in subcategory B, or there is always at least a vertex in subcategory C in all valid hypergraphs. 
\\By the rule of the game, either 1 vertex or 1 hyperedge is removed at a time.

\begin{itemize}
    \item Case 1: Remove a vertex S.
    \\After a vertex S is removed: 
    \begin{itemize}
        \item $|V(CatX)|$ is even.
        \item $|E(CatX)|$ = 1, and $|E(CatY)|$ = 0.
    \end{itemize}
    By definition, the resulting hypergraph $T_{1}$ is an evenly uniform hypergraph with a marked coloring with:
    \begin{itemize}
        \item The resulting vertex set $|V(T_{1})|$ is even.
        \item The resulting hyperedge set $|E(T_{1})|$ is odd.  
    \end{itemize}
    By Theorem 2.8, g($T_{1}$) = 2.

    \item Case 2: Remove a vertex in subcategory A. 
    \\After a vertex in subcategory A is removed: 
    \begin{itemize}
        \item $|V(CatX)|$ is odd, and there is a vertex S (not in category X but in category Y).
        \item $|E(CatX)|$ = 0, and $|E(CatY)|$ is odd.
    \end{itemize}
    By definition, the resulting hypergraph $T_{2}$ is an oddly uniform hypergraph with:
    \begin{itemize}
        \item The resulting vertex set $|V(T_{2})|$ is even.
        \item The resulting hyperedge set $|E(T_{2})|$ is odd.
    \end{itemize}
    By Theorem 2.4, g($T_{2}$) = 3.

    \item Case 3: Remove a vertex in category B (if any). 
    \\After a vertex in category B is removed:
    \begin{itemize}
    \item $|V(CatX)|$ is odd, and there is a vertex S (not in category X but in category Y).
    \item $|E(CatX)|$ = 0, and $|E(CatY)|$ is even.
    \end{itemize}
    By definition, the resulting hypergraph $T_{3}$ is an oddly uniform hypergraph with:
    \begin{itemize}
        \item The resulting vertex set $|V(T_{3})|$ is even.
        \item The resulting hyperedge set $|E(T_{3})|$ is even. 
    \end{itemize}
    By Theorem 2.4, g($T_{3}$) = 0.
    \item Case 4: Remove a vertex in category C (if any). 
    \\After a vertex in category C is removed:
    \begin{itemize}
    \item $|V(CatX)|$ is odd, and there is a vertex S (not in category X but in category Y).
    \item $|E(CatX)|$ = 0, and $|E(CatY)|$ is even.
    \end{itemize}
    
    By definition, the resulting hypergraph $T_{4}$ is an oddly uniform hypergraph with: 
    
    \begin{itemize}
        \item The resulting vertex set $|V(T_{4})|$ is even.
        \item The resulting hyperedge set $|E(T_{4})|$ is even. 
    \end{itemize}
    By Theorem 2.4, g($T_{4}$) = 0.
    
    \item Case 5: Remove a hyperedge in E(CatX).
    \\After a hyperedge in E(CatX) is removed:
    
    \begin{itemize}
    \item $|V(CatX)|$ is even, and there is a vertex S (not in category X but in category Y).
    \item $|E(CatX)|$ = 0, and $|E(CatY)|$ is even.
    \end{itemize}
    By definition, the resulting hypergraph $T_{5}$ is an oddly uniform hypergraph with:
    
    \begin{itemize}
        \item The resulting vertex set $|V(T_{5}|$ is odd.
        \item The resulting hyperedge set $|E(T_{5})|$ is even. 
    \end{itemize}
    By Theorem 2.4, g($T_{5}$) = 1.
    
    \item Case 6: Remove a hyperedge in E(CatY).
    \\After a hyperedge in E(CatY) is removed:
    
    \begin{itemize}
    \item Since $|V(CatX)|$ is even, and there is a vertex S (not in category X but in category Y), $|V(T_{6}|$ is odd.
    \item $|E(CatX)|$ = 1, and $|E(CatY)|$ is odd.
    \end{itemize}
    By Theorem 7, g($T_{6}$) = 1. 

\end{itemize}
Consider all hypergraphs in group III, group IV, and group V. By Sprague-Grundy Theorem, $g(T)$ = mex \{2, 3, 0, 1, 1\} = 4.
\\Thus, consider all hypergraphs in group III, group IV, and group V. By induction, if T is one of the considered hypergraphs where $|E(CatX)|$ = 1 and $|E(CatY)|$ is a positive even number, then g(T) = 4.
\end{paragraph}

\begin{paragraph}{Theorem 9}
\hspace{0.1cm}
\\For all hypergraphs T in group I where $|E(CatX)|$ = 1 and $|E(CatY)|$ is a positive even number, 
\\then g(T) = 0. 

\end{paragraph}

\begin{paragraph}{Proof}
\hspace{1cm}
\\Proof by induction.
\\Base case: Consider a hypergraph T in group I where $|E(CatX)|$ = 1 and $|E(CatY)|$ = 2. 
\begin{itemize}
    \item Since $|V(CatX)|$ is even, and there is a vertex S (not in category X but in category Y). Consequently, the vertex set $|V(T)|$ is odd. 
    \item Since a hypergraph T is in group I, before the game starts, there are only vertices in category A.  
\end{itemize}

By the rule of the game, either 1 vertex or 1 hyperedge is removed
at a time.
\begin{itemize}
    \item Case 1: Remove a vertex S. 
    \\After a vertex S is removed:
    
    \begin{itemize}
        \item $|V(CatX)|$ is even.
        \item $|E(CatX)|$ = 1, and $|E(CatY)|$ = 0.
    \end{itemize}
    
    By definition, the resulting hypergraph $T_{1}$ is an evenly uniform hypergraph with a marked coloring with: 
    
    \begin{itemize}
        \item The resulting vertex set $|V(T_{1})|$ is even.
        \item The resulting hyperedge set $|E(T_{1})|$ is odd. 
    \end{itemize}
    By Theorem 2.8, g($T_{1}$) = 2.

    \item Case 2: Remove a vertex in category A.
    \\After a vertex in category A is removed: 
    
    \begin{itemize}
        \item $|V(CatX)|$ is odd, and there is a vertex S (not in category X but in category Y).
        \item $|E(CatX)|$ = 0, and $|E(CatY)|$ = 1.
    \end{itemize}
    
    By definition, the resulting hypergraph $T_{2}$ is an oddly uniform hypergraph with: 
    
    \begin{itemize}
        \item The resulting vertex set $|V(T_{2})|$ is even.
        \item The resulting hyperedge set $E(T_{2})$ is odd.
    \end{itemize}
    By Theorem 2.4, g($T_{2}$) = 3. 

    \item Case 3: Remove a hyperedge in E(CatX).
    \\After a hyperedge in E(CatX) is removed:
    
    \begin{itemize}
    \item $|V(CatX)|$ is even, and there is a vertex S (not in category X but in category Y).
    \item $|E(CatX)|$ = 0, and $|E(CatY)|$ = 2.
    \end{itemize}
    By definition, the resulting hypergraph $T_{3}$ is an oddly uniform hypergraph with:
    
    \begin{itemize}
        \item The resulting vertex set $|V(T_{3}|$ is odd.
        \item The resulting hyperedge set $|E(T_{3})|$ is even. 
    \end{itemize}
    By Theorem 2.4, g($T_{5}$) = 1.

    \item Case 4: Remove a hyperedge in E(CatY).
    \\After a hyperedge in E(CatY) is removed:
    
    \begin{itemize}
    \item $|V(CatX)|$ is even, and there is a vertex S (not in category X but in category Y).
    \item $|E(CatX)|$ = 1, and $|E(CatY)|$ = 1.
    \end{itemize}
    By Theorem 7, g($T_{4}$) = 1. 
\end{itemize}
By Sprague-Grundy Theorem, $g(T)$ = mex \{2, 3, 1, 1\} = 0.

\vspace{0.1in}

Inductive hypothesis: For all hypergraphs T in group I where $|E(CatX)|$ = 1 and $|E(CatY)|$ = n (n is a positive even number), then g(T) = 0.
\\Consider all hypergraphs T in group I where $|E(CatX)|$ = 1 and $|E(CatY)|$ = n+1 (n+1 is a positive even number). Prove that g(T) = 4.
\begin{itemize}
    \item Since $|V(CatX)|$ is even, and there is a vertex S (not in category X but in category Y). Consequently, the vertex set $|V(T)|$ is odd.
    \item Since a hypergraph T is in group I, before the game starts, there are only vertices in category A. 
\end{itemize}
By the rule of the game, either 1 vertex or 1 hyperedge is removed at a time.

\begin{itemize}
    \item Case 1: Remove a vertex S.
    \\After a vertex S is removed: 
    \begin{itemize}
        \item $|V(CatX)|$ is even.
        \item $|E(CatX)|$ = 1, and $|E(CatY)|$ = 0.
    \end{itemize}
    By definition, the resulting hypergraph $T_{1}$ is an evenly uniform hypergraph with a marked coloring with:
    \begin{itemize}
        \item The resulting vertex set $|V(T_{1})|$ is even.
        \item The resulting hyperedge set $|E(T_{1})|$ is odd.  
    \end{itemize}
    By Theorem 2.8, g($T_{1}$) = 2.

    \item Case 2: Remove a vertex in category A. 
    \\After a vertex in category A is removed: 
    \begin{itemize}
        \item $|V(CatX)|$ is odd, and there is a vertex S (not in category X but in category Y).
        \item $|E(CatX)|$ = 0, and $|E(CatY)|$ is odd.
    \end{itemize}
    By definition, the resulting hypergraph $T_{2}$ is an oddly uniform hypergraph with:
    \begin{itemize}
        \item The resulting vertex set $|V(T_{2})|$ is even.
        \item The resulting hyperedge set $|E(T_{2})|$ is odd.
    \end{itemize}
    By Theorem 2.4, g($T_{2}$) = 3.

    \item Case 3: Remove a hyperedge in E(CatX).
    \\After a hyperedge in E(CatX) is removed:
    
    \begin{itemize}
    \item $|V(CatX)|$ is even, and there is a vertex S (not in category X but in category Y).
    \item $|E(CatX)|$ = 0, and $|E(CatY)|$ is even.
    \end{itemize}
    By definition, the resulting hypergraph $T_{3}$ is an oddly uniform hypergraph with:
    
    \begin{itemize}
        \item The resulting vertex set $|V(T_{3}|$ is odd.
        \item The resulting hyperedge set $|E(T_{3})|$ is even. 
    \end{itemize}
    By Theorem 2.4, g($T_{3}$) = 1.
    
    \item Case 4: Remove a hyperedge in E(CatY).
    \\After a hyperedge in E(CatY) is removed:
    
    \begin{itemize}
    \item Since $|V(CatX)|$ is even, and there is a vertex S (not in category X but in category Y), $|V(T_{6}|$ is odd.
    \item $|E(CatX)|$ = 1, and $|E(CatY)|$ is odd.
    \end{itemize}
    By Theorem 7, g($T_{4}$) = 1. 

\end{itemize}
Consider all hypergraphs T in group I. By Sprague-Grundy Theorem, $g(T)$ = mex \{2, 3, 1, 1\} = 0.
\\Thus, by induction, for all hypergraphs T in group I where $|E(CatX)|$ = 1 and $|E(CatY)|$ is a positive even number, then g(T) = 0. 
\end{paragraph}

\begin{paragraph}{Theorem 10}
\hspace{1cm}
\\For all hypergraphs T in group II where $|E(CatX)|$ = 1 and $|E(CatY)|$ is a positive even number, 
\\then g(T) = 3. 

\end{paragraph}

\begin{paragraph}{Proof}
\hspace{1cm}
\\Proof by induction.
\\Base case: Consider a hypergraph T in group II where $|E(CatX)|$ = 1 and $|E(CatY)|$ = 2. 
\begin{itemize}
    \item Since $|V(CatX)|$ is even, and there is a vertex S (not in category X but in category Y). Consequently, the vertex set $|V(T)|$ is odd. 
    \item Since a hypergraph T is in group II, before the game starts, there are only vertices in category B.  
\end{itemize}

By the rule of the game, either 1 vertex or 1 hyperedge is removed
at a time.
\begin{itemize}
    \item Case 1: Remove a vertex S. 
    \\After a vertex S is removed:
    
    \begin{itemize}
        \item $|V(CatX)|$ is even.
        \item $|E(CatX)|$ = 1, and $|E(CatY)|$ = 0.
    \end{itemize}
    
    By definition, the resulting hypergraph $T_{1}$ is an evenly uniform hypergraph with a marked coloring with: 
    
    \begin{itemize}
        \item The resulting vertex set $|V(T_{1})|$ is even.
        \item The resulting hyperedge set $|E(T_{1})|$ is odd. 
    \end{itemize}
    By Theorem 2.8, g($T_{1}$) = 2.

    \item Case 2: Remove a vertex in category B (if any).
     \\After a vertex in category B is removed:
    \begin{itemize}
    \item $|V(CatX)|$ is odd, and there is a vertex S (not in category X but in category Y).
    \item $|E(CatX)|$ = 0, and $|E(CatY)|$ = 0.
    \end{itemize}
    By definition, the resulting hypergraph $T_{2}$ is an oddly uniform hypergraph with:
    \begin{itemize}
        \item The resulting vertex set $|V(T_{2})|$ is even.
        \item The resulting hyperedge set $|E(T_{2})|$ is even. 
    \end{itemize}
    By Theorem 2.4, g($T_{2}$) = 0.
    
    \item Case 3: Remove a hyperedge in E(CatX).
    \\After a hyperedge in E(CatX) is removed:
    
    \begin{itemize}
    \item $|V(CatX)|$ is even, and there is a vertex S (not in category X but in category Y).
    \item $|E(CatX)|$ = 0, and $|E(CatY)|$ = 2.
    \end{itemize}
    By definition, the resulting hypergraph $T_{3}$ is an oddly uniform hypergraph with:
    
    \begin{itemize}
        \item The resulting vertex set $|V(T_{3}|$ is odd.
        \item The resulting hyperedge set $|E(T_{3})|$ is even. 
    \end{itemize}
    By Theorem 2.4, g($T_{3}$) = 1.

    \item Case 4: Remove a hyperedge in E(CatY).
    \\After a hyperedge in E(CatY) is removed:
    
    \begin{itemize}
    \item $|V(CatX)|$ is even, and there is a vertex S (not in category X but in category Y).
    \item $|E(CatX)|$ = 1, and $|E(CatY)|$ = 1.
    \end{itemize}
    By Theorem 7, g($T_{4}$) = 1. 
\end{itemize}
By Sprague-Grundy Theorem, $g(T)$ = mex \{2, 0, 1, 1\} = 3.

\vspace{0.1in}

Inductive hypothesis: For all hypergraphs T in group II where $|E(CatX)|$ = 1 and $|E(CatY)|$ = n (n is a positive even number), then g(T) = 4.
\\Consider all hypergraphs T in group II where $|E(CatX)|$ = 1 and $|E(CatY)|$ = n+1 (n+1 is a positive even number). Prove that g(T) = 3.
\begin{itemize}
    \item Since $|V(CatX)|$ is even, and there is a vertex S (not in category X but in category Y). Consequently, the vertex set $|V(T)|$ is odd.
    \item Since a hypergraph T is in group II, before the game starts, there are only vertices in category B.
\end{itemize}
By the rule of the game, either 1 vertex or 1 hyperedge is removed at a time.

\begin{itemize}
    \item Case 1: Remove a vertex S.
    \\After a vertex S is removed: 
    \begin{itemize}
        \item $|V(CatX)|$ is even.
        \item $|E(CatX)|$ = 1, and $|E(CatY)|$ = 0.
    \end{itemize}
    By definition, the resulting hypergraph $T_{1}$ is an evenly uniform hypergraph with a marked coloring with:
    \begin{itemize}
        \item The resulting vertex set $|V(T_{1})|$ is even.
        \item The resulting hyperedge set $|E(T_{1})|$ is odd.  
    \end{itemize}
    By Theorem 2.8, g($T_{1}$) = 2.

    \item Case 2: Remove a vertex in category B. 
    \\After a vertex in category B is removed:
    \begin{itemize}
    \item $|V(CatX)|$ is odd, and there is a vertex S (not in category X but in category Y).
    \item $|E(CatX)|$ = 0, and $|E(CatY)|$ is even.
    \end{itemize}
    By definition, the resulting hypergraph $T_{2}$ is an oddly uniform hypergraph with:
    \begin{itemize}
        \item The resulting vertex set $|V(T_{2})|$ is even.
        \item The resulting hyperedge set $|E(T_{2})|$ is even. 
    \end{itemize}
    By Theorem 2.4, g($T_{2}$) = 0.

    \item Case 3: Remove a hyperedge in E(CatX).
    \\After a hyperedge in E(CatX) is removed:
    
    \begin{itemize}
    \item $|V(CatX)|$ is even, and there is a vertex S (not in category X but in category Y).
    \item $|E(CatX)|$ = 0, and $|E(CatY)|$ is even.
    \end{itemize}
    By definition, the resulting hypergraph $T_{3}$ is an oddly uniform hypergraph with:
    
    \begin{itemize}
        \item The resulting vertex set $|V(T_{3}|$ is odd.
        \item The resulting hyperedge set $|E(T_{3})|$ is even. 
    \end{itemize}
    By Theorem 2.4, g($T_{3}$) = 1.
    
    \item Case 4: Remove a hyperedge in E(CatY).
    \\After a hyperedge in E(CatY) is removed:
    
    \begin{itemize}
    \item Since $|V(CatX)|$ is even, and there is a vertex S (not in category X but in category Y), $|V(T_{4}|$ is odd.
    \item $|E(CatX)|$ = 1, and $|E(CatY)|$ is odd.
    \end{itemize}
    By Theorem 7, g($T_{4}$) = 1. 

\end{itemize}
For all hypergraphs in group II, by Sprague-Grundy Theorem, $g(T)$ = mex \{2, 0, 1, 1\} = 3.
\\Thus, by induction, for all hypergraphs in group II where $|E(CatX)|$ = 1 and $|E(CatY)|$ is a positive even number, then g(T) = 3.
\end{paragraph}
\begin{paragraph}{}
\paragraph{Part 3: Applications.}
\begin{itemize}
    \item A game on any geometric object whose Planar Graph satisfying all of the requirements in part 1 can be played by using the winning strategy in part 2.
    \\Example: Graphs.
    \item A game on any geometric object whose Schlegel diagram satisfying all of the requirements in part 1 can be played by using the winning strategy in part 2.
    \\Example: Pyramids.
    \item A game on any geometric object whose graphical projection satisfying all of the requirements in part 1 can be played by using the winning strategy in part 2.
    \\Example: Chocolate boxes.
    \item A game on any geometric object whose n-polytope satisfying all of the requirements in part 1 can be played by using the winning strategy in part 2.
    \\Example: polygons that are self-intersecting with varying densities of different regions.
    \item A game on any complex polygons satisfying all of the requirements in part 1 can be played by using the winning strategy in part 2.
    \item A game on any polyhedrons which have the net of the polyhedrons satisfying all of the requirements in part 1 can be played by using the winning strategy in part 2.
\end{itemize}
\end{paragraph}

\section*{Acknowledgments}
This article was served as T. H. Molena's Graduating Senior Capstone from Berea College, which is a liberal arts college in Berea, Kentucky, in August of 2021.

\bibliographystyle{unsrt}  
\bibliography{references}  

\end{document}